\documentclass[leqno]{amsart}

\usepackage[foot]{amsaddr}

\usepackage{color}
\usepackage{geometry}
\usepackage{subfiles}
\usepackage{graphicx}
\usepackage[section]{placeins}
\usepackage{hyperref}
\usepackage{stmaryrd}
\usepackage{amsmath,amssymb,amsfonts,amsthm,textcomp}
\usepackage{mathtools}
\usepackage{tensor}
\usepackage{multicol}
\usepackage{subfig}
\graphicspath{
              {./figures/}
              {./figures_section-1/}
              {./figures_section-2/}
              {./figures_section-3/}
             }

\newfont{\tenbfsl}{cmbxti9 scaled 1200}
\newfont{\tenbbb}{msbm10}
\newfont{\svnbbb}{msbm8}

\newcommand{\bs}[1]{\boldsymbol{#1}}

\newcommand{\cl}[1]{\mathcal{#1}}

\newcommand{\bb}[1]{\mathbb{#1}}

\newcommand{\norm}[2]{\left\lVert {#1} \right\rVert_{#2}}

\newcommand{\Grad}{\nabla}

\theoremstyle{remark}

\theoremstyle{definition}

\newcounter{syn}[section] 
\renewcommand{\thesyn}{\arabic{section}.\arabic{syn}}

\makeatletter
\def\threevdots{\mskip+4mu\vbox{\baselineskip2.25\p@ \lineskiplimit\z@
  \kern4.9\p@\hbox{.}\hbox{.}\hbox{.}}\mskip+3.8mu}
\makeatother

\begin{document}

\title{Deep NURBS---Admissible Physics-informed Neural Networks}
\author{Hamed Saidaoui$^{\sharp}$, Luis Espath$^1$ \& Ra\'ul Tempone$^{2,3,4}$}
\address{$^1$School of Mathematical Sciences, University of Nottingham, Nottingham, NG7 2RD, United Kingdom}
\address{$^2$Department of Mathematics, RWTH Aachen University, Geb\"{a}ude-1953 1.OG, Pontdriesch 14-16, 161, 52062 Aachen, Germany.}
\address{$^3$King Abdullah University of Science \& Technology (KAUST), Computer, Electrical and Mathematical Sciences \& Engineering Division (CEMSE), Thuwal 23955-6900, Saudi Arabia.}
\address{$^4$Alexander von Humboldt Professor in Mathematics for Uncertainty Quantification, RWTH Aachen University, Germany.}
\email{$^\sharp$ hamed.saidaoui@kaust.edu.sa}

\date{\today}

\begin{abstract}
\noindent
In this study, we propose a new numerical scheme for physics-informed neural networks (PINNs) that enables precise and inexpensive solutions for partial differential equations (PDEs) in case of arbitrary geometries while strongly enforcing Dirichlet boundary conditions. The proposed approach combines admissible NURBS parametrizations (admissible in the calculus of variations sense, that is, satisfying the boundary conditions) required to define the physical domain and the Dirichlet boundary conditions with a PINN solver. Therefore, the boundary conditions are automatically satisfied in this novel Deep NURBS framework. Furthermore, our sampling is carried out in the parametric space and mapped to the physical domain. This parametric sampling works as an importance sampling scheme since there is a concentration of points in regions where the geometry is more complex. We verified our new approach using two-dimensional elliptic PDEs when considering arbitrary geometries, including non-Lipschitz domains. Compared to the classical PINN solver, the Deep NURBS estimator has a remarkably high accuracy for all the studied problems. Moreover, a desirable accuracy was obtained for most of the studied PDEs using only one hidden layer of neural networks. This novel approach is considered to pave the way for more effective solutions for high-dimensional problems by allowing for a more realistic physics-informed statistical learning framework to solve PDEs. 
\\
\textbf{AMS subject classifications:}
$\cdot$
35L65 
$\cdot$

\end{abstract}

\maketitle

\tableofcontents                        


\section{Introduction}

Deep Learning (DL) exhibits unprecedented advancement in the last two decades \cite{LeCun2015}. Its usage in numerous disciplines has resulted in diverse successful implementations such as language processing \cite{speech2012, speech2013} and image recognition \cite{img_recg2012,img_recg2015}. With the exception of a few studies, the application of machine learning (ML) along with neural networks (NNs) in general, in particular to scientific problems, was not as structured as in the previous topics (for example \cite{lagaris1998, carleo2017,carleo2018, Torlai2018, schmidt2009, snyder2012}). 

Among these few efforts was the seminal work of Lagaris \cite{lagaris1998}, who applied very simplistic NN models to solve several ordinary and partial differential equations (PDEs) of different orders (the accuracy was order-dependent though). This concept has been revived very recently by many researchers \cite{raissi2019,berg2018,weinan2018,sirignano2018} who used a very similar concept and coined it as "physics-informed neural networks (PINNs)" (two such examples were deep Ritz and deep Galerkin for the last two references, respectively). Since then, PINNs have gained considerable fame because of their simple implementation and their concept that allows for a combination of NNs and the already established theories. PINNs have enabled the shift to the unsupervised ML scheme owing to the laws and constraints implemented within its framework. Targeting the solution of PDEs, PINNs have been proven to be very efficient for many complicated and challenging problems (\cite{cai2021, Kurth2018}) and for different types of PDEs. Furthermore, the convenience of the accompanying sampling methods made it adequate for problems with high-dimensional domains. Because PINNs have been specifically proposed for methods that consider the collocation points as training data with a discrete loss term; note that other variants relying on the same principle coexisted, e.g., deep Ritz \cite{weinan2018, samaniego2020} and deep Galerkin \cite{sirignano2018}. The first method uses the weak form of the PDE in which the variational problem (energy minimization) can be solved using stochastic optimizers and Monte Carlo sampling techniques. In the deep Galerkin scheme, the gradients have been taken care of using a tailored stochastic approach. 

Many variations have been proposed and proven to be effective in specific cases, in addition to these three primary methods. Inverse problems have been tackled using stochastic differential equations (SDEs) in this reference \cite{Han2018}. Lu et al. \cite{Lu2021} proposed a method for learning non-linear operators, while Cai et al. \cite{cai2021} have implemented a modified version of PINNs to infer electroconvection in multiphysics systems. A comprehensive review on ML for fluid mechanics is presented in \cite{Bru20}. 
The mathematical formulation of PINNs has been discussed in certain relevant studies from an uncertainty quantification perspective. Mishra et al. \cite{mishra2022} provided an estimate for its generalization error. They primarily provided bounds for PINNs based on quadrature and random sampling based-PINNs. Fang et al. \cite{fang2021} established convergence rates for the NNs-based PDE solver. To develop an approximation for the derivatives, they combined the NNs and the differential operator (as it is used in Finite element methods).

As a trend in applying ML tools to scientific and engineering problems, we naturally wonder what makes PINNs and their variants (deep Ritz and deep Galerkin) as efficient and vulnerable as they are (or considered to be). To address this question, one should examine what makes these methods different from answering this question. Surely its simple concept makes it easy to implement, particularly for complicated tasks in which simplicity plays an important role in its success. One has to admit that the nature of its loss function, which incorporates the system's physical laws, has to do with its convergence rate to the sought solution. The auto-differentiation, being an exact, differentiating tool, is, in turn, at the heart of PINNs success. 

However, PINNs have observed many limitations when it comes to problems with discontinuities and computationally demanding problems \cite{weinan2018}. These problems are very demanding in terms of NN depth and computational work; the latter two aspects make the PINNs usage challenging. From this perspective, one might question what has to change to allow PINNs to cope with discontinuities without being extremely expensive. Is it the way the loss function is presented (weak form \cite{weinan2018, samaniego2020} or strong form \cite{raissi2019}), or is it related to the way the gradients are calculated (entirely auto-differentiation - based or hybrid \cite{fang2021}). If we examine the previous references, one can conclude that none of the previous aspects is important in mitigating the complexity of the NN architecture. 

\begin{figure}
   \centering
   \subfloat[Admissible (ansatz) NURBS $\varphi(\bs{x})$]{\includegraphics[width=0.33\textwidth]{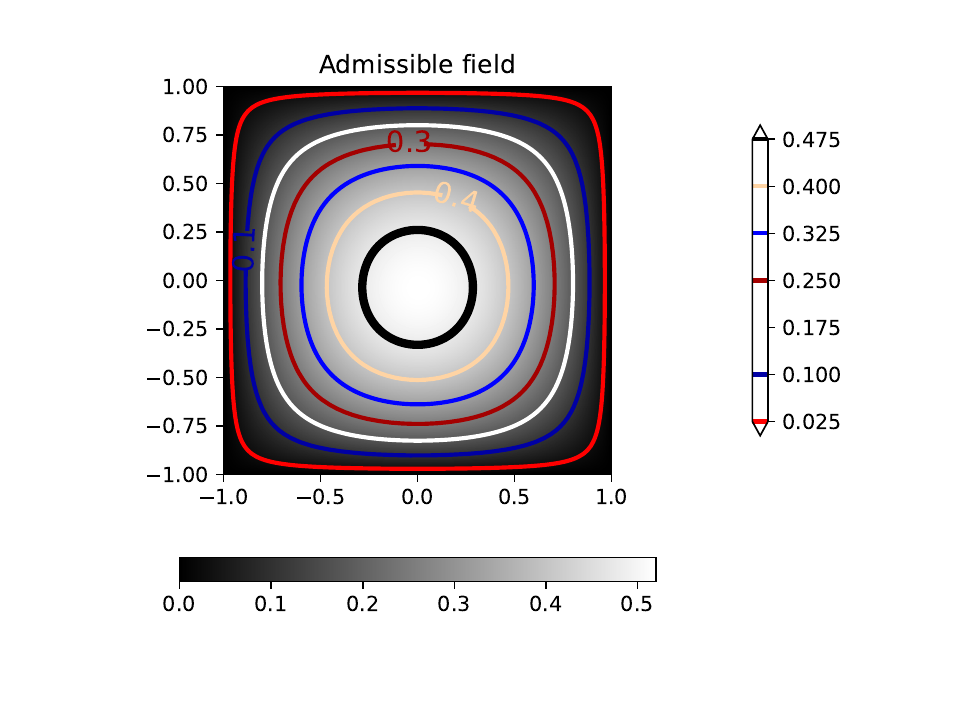}}
   \subfloat[Neural Network $NN(\bs{x})$]{\includegraphics[width=0.33\textwidth]{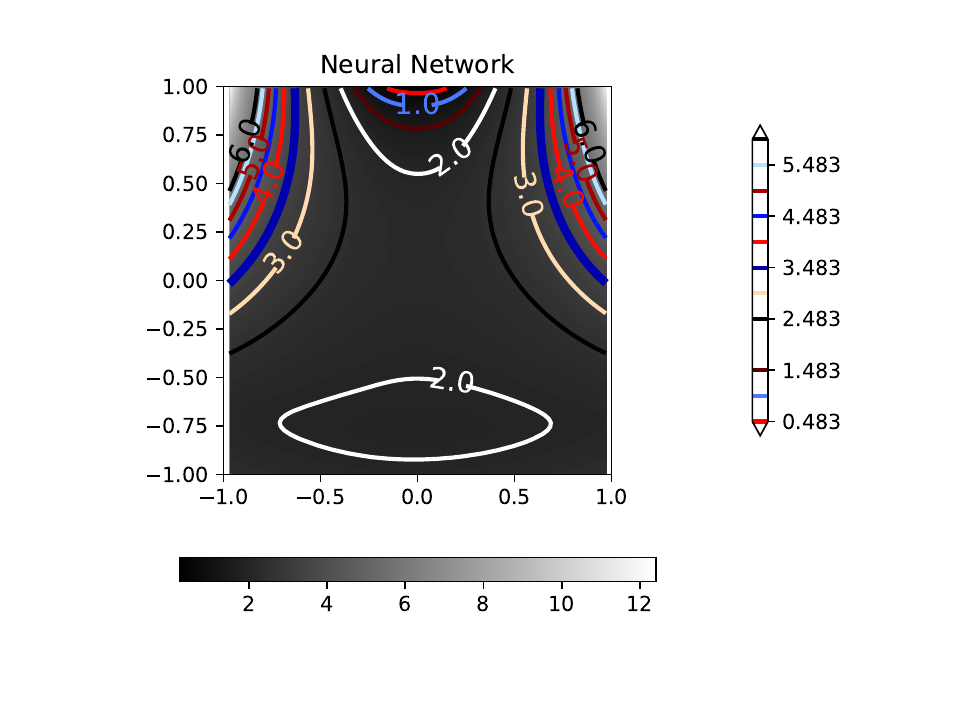}}
   \subfloat[Admissible Neural Network $NN(\bs{x})\varphi(\bs{x})$]{\includegraphics[width=0.33\textwidth]{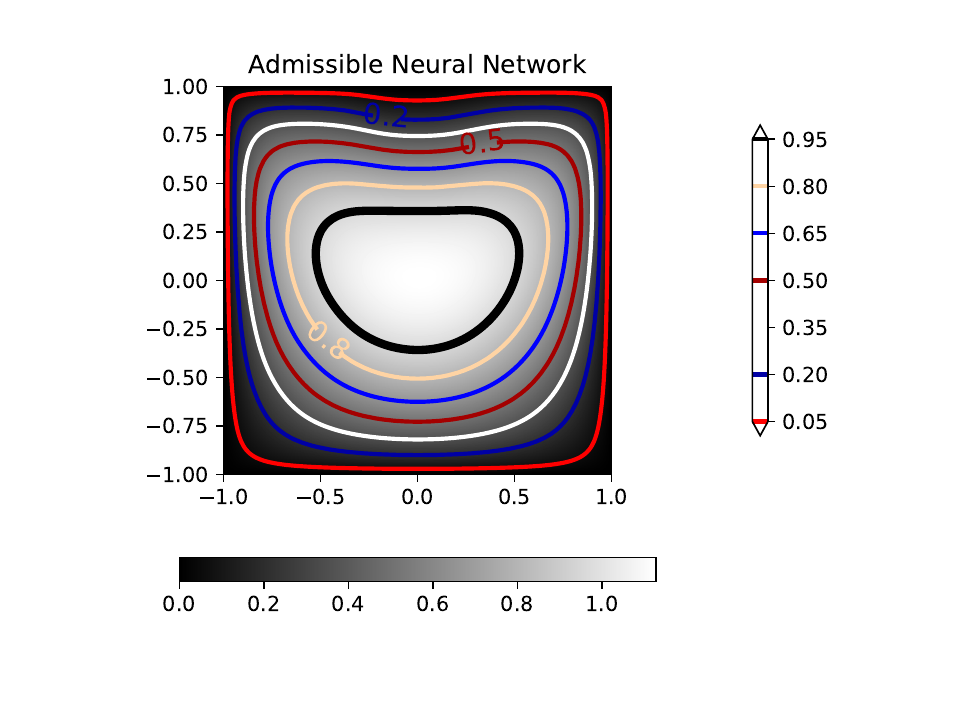}}
   \caption{Exemplification of Deep NURBS' method. From left to right: Admissible NURBS (ansatz) $\varphi(\bs{x})$, Neural Network $NN(\bs{x})$, and the product of the two first, that is, Admissible Neural Network $NN(\bs{x})\varphi(\bs{x})$, for the homogeneous Dirichlet case.}
   \label{fg:example}
\end{figure}

Shin et al. \cite{shin2020} demonstrated the sequence of minimizers generated in the stochastic optimization of PINNs strongly converges to the solution of the PDE in $C^0$. However, if the boundary conditions are satisfied, the minimizers converge to the sought solution in $H^1$ instead of $C^0$. This insight has been supported by recent works \cite{sukumar2022, berg2018} and genuinely by the original work of Lagaris \cite{lagaris1998}. The Lagaris approach was to use functions that naturally fulfill the boundary conditions. Although this approach proved to be computationally very efficient, in terms of accuracy, this approach is not practical regarding complex geometries with non-homogeneous boundaries. 
From this last point, we follow the concept described in Lagaris' paper. Nevertheless, we use admissible (in the sense of calculus of variations, that is, satisfying the boundary conditions) non-uniform rational B-splines (NURBS) parameterizations, which are similarly used in IsoGeometric Analysis, as an efficient approach to impose Dirichlet boundary conditions. To our knowledge, our work is the first to impose boundary conditions using admissible NURBS in the PINNs framework. We call our method Deep-NURBS, which is a combination of deep neural networks (DNNs) and NURBS. We use the weak form of the loss function (energy term), and we implement our algorithms in Tensorflow \cite{tensorflow2016}. Figure \ref{fg:example} shows the admissible NURBS (which works as an ansatz), the NN, and the product of the two, rendering the Admissible NN that ultimately leads to the Deep-NURBS method. Moreover, despite being an ML-based method, we demonstrate that our method is robust against slight changes in the NN parameters. This fact allows us to overcome computational costs coming with hyperparameter optimization. Lastly, note that the NURBS parameterization does not change the overall computational cost, that is, our Deep NURBS method increases the approximability of the Neural Networks without increasing the overall cost for a fix Neural Network.

In this study, we overview NURBS parameterizations in Section \ref{sec:sec1}. In Section \ref{sec:sec2}, we discuss our approach based on combining DNNs with NURBS. Section \ref{sec:sec3} will be devoted to discussing our result of applying deep NURBS to various problems with non-trivial complex geometries. Finally, we wrap up this study with a conclusion.

\section{Mathematical background}\label{sec:sec1}

Our notation is as follows: scalar and vectorial fields are denoted by lower-case letters and lower-case boldface letters, respectively. The gradient and Laplacian operators are denoted by $\Grad$ and $\triangle$, respectively. The physical domain is denoted as $\cl{D}$ with boundary $\partial\cl{D}$. Moreover, we let $\partial\cl{D}\coloneqq\partial\cl{D}_{\mathrm{D}}\cup\partial\cl{D}_{\mathrm{N}}$ where the Dirichlet boundary conditions are considered on $\cl{D}_{\mathrm{D}}$ while the Neumann type of boundary conditions are considered on $\cl{D}_{\mathrm{N}}$. We construct our approximation on the Sobolev space of all square-integrable functions on $\cl{D}$ with square-integrable derivatives. Thus, we denote the Lebesgue space as $L^2$ and $\bs{L}^2$ for scalar-valued and vector-valued functions, respectively. Similarly, we denote the Sobolev space of scalar- and vector-valued functions as $H^k$ and $\bs{H}^k$, respectively, where the norm of the $k$th gradient is $\bs{L}^2$.

Let us consider two cases. Let $\cl{L}[\bs{u}(\bs{x};\bs{\theta})]$ be (a) the pointwise residual of a PDE or (b) the energy density of the underlying system given an ansatz $\bs{u}$ parametrized by $\bs{\theta} \in \bb{R}^Q$ where $Q$ is the total number of degrees of freedom (e.g., the total number of neural network parameters). 
Next, let $\bs{u}^\ast$ be the solution of this PDE with values $\bs{u}_\mathrm{D}$, on $\partial\cl{D}_{\mathrm{D}}$ while its normal derivative $\partial_n\bs{u}^\ast$ takes the value $\bs{g}_{\mathrm{N}}$ on $\partial\cl{D}_{\mathrm{N}}$. In a learning type of framework, we may state problems (a) and (b) as follows,
\begin{equation}\label{eq:usual.pinn}
\bs{\theta}^\ast\coloneqq\underset{\bs{\theta}\in \bb{R}^Q}{\arg\min} \, \norm{\cl{L}[\bs{u}(\cdot;\bs{\theta})]}{L^2(\cl{D})}^2 + \mathrm{Dirichlet\ terms\ on\ } \partial\cl{D}_{\mathrm{D}} + \mathrm{Neumann\ terms\ on\ } \partial\cl{D}_{\mathrm{N}}.
\end{equation}

Note that the Dirichlet term is imposed in a weak sense where authors in Ref. \cite{weinan2018} use a penalization term to impose the Dirichlet boundary conditions.
Instead, our approach is to impose the Dirichlet boundary conditions in a strong sense. To this end, we consider a parameterization of our domain $\cl{D}$,
\begin{equation}
\bs{x}\coloneqq\bs{\chi}(\bs{\xi}),
\end{equation}
where we assume that $\bs{\chi}$ is a bijective mapping with $\bs{x}\in\cl{D}$ and $\bs{\xi}$ is a parameter that lives in the parameter space $\bs{\Xi}$. We here use the notion of admissible fields. A family of functions is said to be admissible if Dirichlet boundary conditions are satisfied. Consider the following family of admissible vector-valued parameterizations
\begin{equation}
\cl{Z}\coloneqq\{\bs{\zeta}\,|\,\bs{\zeta}\in\bs{H}^k(\cl{D})\quad\wedge\quad\bs{\zeta}(\bs{x})=\bs{u}_\mathrm{D}(\bs{x}),\quad\forall\,\bs{x}=\bs{\chi}(\bs{\xi})\in\partial\cl{D}_{\mathrm{D}}\}.
\end{equation}
Moreover, we define a family of auxiliary smooth scalar-valued functions $\cl{P}$ with vanishing boundary, i.e., $\varphi\colon\bs{\chi}(\bs{\Xi})\mapsto\bb{R}$ with $\varphi(\bs{\chi}(\bs{\xi}))=0$ for all $\bs{\xi}=\bs{\chi}^{-1}(\bs{x})$ such that $\bs{x}=\bs{\chi}(\bs{\xi})\in\partial\cl{D}_{\mathrm{D}}$. Furthermore, we require $\varphi\in{H}^k(\cl{D})$. Thus,
\begin{equation}
\cl{P}\coloneqq\{\varphi\,|\,\varphi\in{H}^k(\cl{D})\quad\wedge\quad\varphi(\bs{x})=0,\quad\forall\,\bs{x}=\bs{\chi}(\bs{\xi})\in\partial\cl{D}_{\mathrm{D}}\}.
\end{equation}

The problem statement \eqref{eq:usual.pinn} reads
\begin{equation}\label{eq:admissible.NURBS.pinn}
\left\{
\begin{aligned}
&\bs{\theta}^\ast\coloneqq\underset{\bs{\theta}\in \bb{R}^Q}{\arg\min} \, \norm{\cl{L}[\varphi(\cdot)\mskip2mu \bs{u}(\cdot;\bs{\theta})+\bs{\zeta}(\cdot)]}{L^2(\cl{D})}^2 + \mathrm{Neumann\ terms\ on\ } \partial\cl{D}_{\mathrm{N}},\\[4pt]
&\text{given an arbitrary }\varphi\in\cl{P}\text{ and }\bs{\zeta}\in\cl{Z}.
\end{aligned}
\right.
\end{equation}
Furthermore, for the scalar case, that is, we replace the vector-valued function $\bs{u}$ with the scalar valued function $u$, expression \eqref{eq:admissible.NURBS.pinn} takes the following form
\begin{equation}\label{eq:admissible.NURBS.pinn.scalar.nonhomogeneous}
\left\{
\begin{aligned}
&\bs{\theta}^\ast\coloneqq\underset{\bs{\theta}\in \bb{R}^Q}{\arg\min} \, \norm{\cl{L}[\varphi(\cdot)\mskip2muu(\cdot;\bs{\theta})+\zeta(\cdot)])}{L^2(\cl{D})}^2 + \mathrm{Neumann\ terms\ on\ } \partial\cl{D}_{\mathrm{N}},\\[4pt]
&\text{given an arbitrary }\varphi\in\cl{P}\text{ and }\zeta\in\cl{Z},
\end{aligned}
\right.
\end{equation}
where here $\cl{Z}$ is
\begin{equation}
\cl{Z}\coloneqq\{\zeta\,|\,\zeta\in H^k(\cl{D})\quad\wedge\quad\zeta(\bs{x})=u_\mathrm{D}(\bs{x}),\quad\forall\,\bs{x}=\bs{\chi}(\bs{\xi})\in\partial\cl{D}_{\mathrm{D}}\}.
\end{equation}
If we instead restrict attention to the scalar case with homogeneous boundary conditions, the function $\varphi$ becomes our admissible field such that the expression \eqref{eq:admissible.NURBS.pinn.scalar.nonhomogeneous} specializes to
\begin{equation}\label{eq:admissible.NURBS.pinn.scalar.homogeneous}
\left\{
\begin{aligned}
&\bs{\theta}^\ast\coloneqq\underset{\bs{\theta}\in \bb{R}^Q}{\arg\min} \, \norm{\cl{L}[\varphi(\cdot)\mskip2muu(\cdot;\bs{\theta})])}{L^2(\cl{D})}^2 + \mathrm{Neumann\ terms\ on\ } \partial\cl{D}_{\mathrm{N}},\\[4pt]
&\text{given an arbitrary }\varphi\in\cl{P}.
\end{aligned}
\right.
\end{equation}

We will demonstrate how $\varphi$ and $\zeta$ may be developed using NURBS parameterizations.

\subsection{Construction of admissible NURBS parameterizations}

Now, let us consider the subset of NURBS functions $\cl{P}^\lambda\subset\cl{P}$ and $\cl{Z}^\lambda\subset\cl{Z}$ which will work as an ansatz. The Cox--deBoor recursive formulation \cite{Deb72,Cox72} is usually adopted to evaluate B-\emph{spline} basis functions based on \cite{Bez72}. In $d$ dimensions, B-\emph{splines} are obtained considering a given knot vector $\bs{\Xi}=\bigotimes_{\jmath=1}^d\Xi^\jmath$ defined over the parametric space with polynomial degree $p_\jmath$ along the parametric direction $\xi^\jmath$. Finally, to map from the parametric space to the physical one, we use $n_\jmath+1$ control points. The NURBS basis functions and physical domain are discussed below:
\begin{itemize}
    \item[Bases:] B-\emph{splines} in the $\jmath$th spatial direction
    \begin{equation}\label{eq:N}
    \begin{split}
    & N_{i, 0} \left( \xi \right) = \begin{cases}
    1& \text{if $\xi_i \leqslant \xi < \xi_{i+1}$},\\
    0& \text{otherwise}.
    \end{cases}\\
    & N_{i, p} \left( \xi \right) = \dfrac{ \xi - \xi_i }{ \xi_{i+p} - \xi_i} N_{i, p-1} \left( \xi \right)+ \dfrac{ \xi_{i+p+1} - \xi }{ \xi_{i+p+1} - \xi_{i+1}} N_{i+1, p-1} \left( \xi \right)
    \end{split}
    \end{equation}
    with
    \begin{equation}\label{eq:knots}
    \Xi^\jmath =\{\underbrace{0, \ldots, 0}_{p_\jmath+1}, \xi^\jmath_{p_\jmath+1}, \ldots, \xi^\jmath_{s_\jmath-p_\jmath-1}, \underbrace{1, \ldots, 1}_{p_\jmath+1}\},
    \end{equation}
    where $s_\jmath = n_\jmath + p_\jmath + 1$. $n_\jmath$ is the number of basis functions and $p_\jmath$ its degree along direction $\jmath$.
    \item[Domain:]
    Given $n_\jmath+1$ basis functions $N_{i,p_\jmath}(\xi^\jmath)$ of degree $p_\jmath$ along the $\jmath$th direction and the corresponding control points $\bs{p}$, a NURBS curve, surface, and volume may be described as follows:
    \begin{equation}
    \bs{x}(\bs{\xi})=\dfrac{1}{W(\bs{\xi})}\sum_{i=0}^{n_1} N_{i,p_1}(\xi^1)\bs{p}_{i},
    \end{equation}
    \begin{equation}\label{eq:NURBS2d}
    \bs{x}(\bs{\xi})=\dfrac{1}{W(\bs{\xi})}\sum_{i=0}^{n_1}\sum_{j=0}^{n_2} N_{i,p_1}(\xi^1)N_{j,p_2}(\xi^2)\bs{p}_{i,j},
    \end{equation}
    and
    \begin{equation}
    \bs{x}(\bs{\xi})=\dfrac{1}{W(\bs{\xi})}\sum_{i=0}^{n_1}\sum_{j=0}^{n_2}\sum_{k=0}^{n_3} N_{i,p_1}(\xi^1)N_{j,p_2}(\xi^2)N_{k,p_3}(\xi^3)\bs{p}_{i,j,k},
    \end{equation}
    respectively. The weights $W$ are respectively defined:
    \begin{equation}
    W(\bs{\xi})=\sum_{i=0}^{n_1} N_{i,p_1}(\xi^1),
    \end{equation}
    \begin{equation}
    W(\bs{\xi})=\sum_{i=0}^{n_1}\sum_{j=0}^{n_2} N_{i,p_1}(\xi^1)N_{j,p_2}(\xi^2),
    \end{equation}
    and
    \begin{equation}
    W(\bs{\xi})=\sum_{i=0}^{n_1}\sum_{j=0}^{n_2}\sum_{k=0}^{n_3} N_{i,p_1}(\xi^1)N_{j,p_2}(\xi^2)N_{k,p_3}(\xi^3).
    \end{equation}
\end{itemize}

Our ansatz for the admissible construction follows the same reasoning. For example, in two dimensions, scalar fields are given as follows: 
\begin{equation}\label{eq:admField}
\varphi(\bs{\xi})=\dfrac{1}{W(\bs{\xi})}\sum_{i=0}^{n_1}\sum_{j=0}^{n_2} N_{i,p_1}(\xi^1)N_{j,p_2}(\xi^2)\overline{\varphi}_{i,j},
\end{equation}
with scalar coefficients $\overline{\varphi}_{i,j}$, while vector fields
\begin{equation}
\bs{\zeta}(\bs{\xi})=\dfrac{1}{W(\bs{\xi})}\sum_{i=0}^{n_1}\sum_{j=0}^{n_2} N_{i,p_1}(\xi^1)N_{j,p_2}(\xi^2)\overline{\bs{\zeta}}_{i,j},
\end{equation}
are expressed with vector coefficients $\overline{\bs{\zeta}}_{i,j}$.

Finally, for the gradient of the admissible fields with $\bs{\xi}=\bs{\chi}^{-1}(\bs{x})$, we have
\begin{equation}
\Grad_{\bs{\xi}}\varphi(\bs{\xi})\coloneqq(\Grad_{\bs{\xi}}\bs{x})\Grad_{\bs{x}}\varphi(\bs{\chi}^{-1}(\bs{x})).
\end{equation}
Thus,
\begin{equation}
\Grad_{\bs{x}}\varphi(\bs{\chi}^{-1}(\bs{x}))=(\Grad_{\bs{\xi}}\bs{x})^{-1}\Grad_{\bs{\xi}}\varphi(\bs{\xi}),
\end{equation}
for scalar fields. Similarly,
\begin{equation}
\Grad_{\bs{x}}\bs{\zeta}(\bs{\chi}^{-1}(\bs{x}))=(\Grad_{\bs{\xi}}\bs{x})^{-1}\Grad_{\bs{\xi}}\bs{\zeta}(\bs{\xi}),
\end{equation}
for vector fields.

\section{Deep NURBS}\label{sec:sec2}

In this section, we will show the foundations of Deep NURBS. We start by stating the partial differential equation to be solved (in its abstract form). We then combine the NURBS admissible field parametrizations that we have encountered in the previous section with the neural network framework of the PINNs. 

\subsection{Partial differential equation and method overview}

Let us consider a physical domain $\cl{D}$ embedeed in $\mathbb{R}^d$ and let $\partial \cl{D}$ denote the boundaries limiting the domain $\cl{D}$. The general form of a partial differential equation defined on the domain $\cl{D}$ with homogeneous Dirichlet boundary conditions reads as follows:
\begin{equation}\label{eq:pde}
\begin{aligned}
\cl{L}[u](x) &= f \quad \bs{x} \in \cl{D} , \\
B[u](x) &= 0   \quad \bs{x} \in \partial \cl{D},
\end{aligned}
\end{equation}
$\cl{L}$ and $B$ are the differential and boundary operators, respectively, and $u$ is the sought solution. The PDE above can be turned into a weak form.
For example, one can consider the two-dimensional Poisson equation defined over the domain $\cl{D}$: 
\begin{equation}\label{eq:poisson1}
\begin{aligned}
-\Delta u(\bs{x})=f(\bs{x}), & \quad \bs{x} \in \cl{D}, \\
        u(\bs{x})=0, & \quad \bs{x} \in \partial \cl{D},
\end{aligned}
\end{equation}
where its weak form leads to the minimization of the energy term defined as \cite{weinan2018}:
\begin{equation}\label{eq:weak_form}
\min _{u \in H^{1}} \int_{\cl{D}}\left(\frac{1}{2}\left|\nabla_{\bs{x}} u(\bs{x})\right|^{2}-f(\bs{x}) u(\bs{x})\right) d \bs{x}+\lambda \int_{\partial \cl{D}}u^{2}(s) d s, 
\end{equation}
Here, $H^1$ is the Sobolev space for functions with square integrable derivatives and $\lambda$ is a penalizing parameter to account for the boundary conditions expressed by the second integral of the equation. The sought solution of the PDE $u(\bs{x})$ is given as a combination between a NURBS and a neural networks function. We will provide more insights on how to construct $u$ in what follows. Note that in the current study, we will only consider the first integral of \eqref{eq:weak_form} as our approximation is an admissible Neural Network, that is, our approximation satisfy the Dirichlet boundary conditions by construction. Also, note that the first integral of Eq. \eqref{eq:weak_form} can be seen as an expectation of the integrand with respect to a uniform density. The minimization reads then as follows :
\begin{equation}
\min _{u \in H^{1}, \bs{x}\sim \mathbb{U}([0,1]^2)} \mathbb{E}_{\bs{x}}\left[\frac{1}{2}\left|\nabla_{\bs{x}} u(\bs{x})\right|^{2}-f(\bs{x}) u(\bs{x})\right],
\end{equation}

\noindent In practice, a more careful sampling scheme should be followed for better results. In section \ref{sec:sampl_mthd}, we will see that NURBS geometry generation allows for a spontaneous importance sampling scheme. Compared to classical machine learning solvers, PINNs (and their variants) categorize under the umbrella of unsupervised methods because no labeled data is required to solve the PDE. The only exception is the boundary points which are paramount for the problem's well-posedness. However, the method we propose in the current manuscript does not even need boundary information to solve the PDE at hand, and thus it can be classified as a fully unsupervised scheme. Indeed one pre-imposes hard boundary conditions from the beginning. 

\subsection{Physics - Informed Neural Networks (PINNs)}

\begin{figure}
    \centering
    \includegraphics[width=0.7\textwidth]{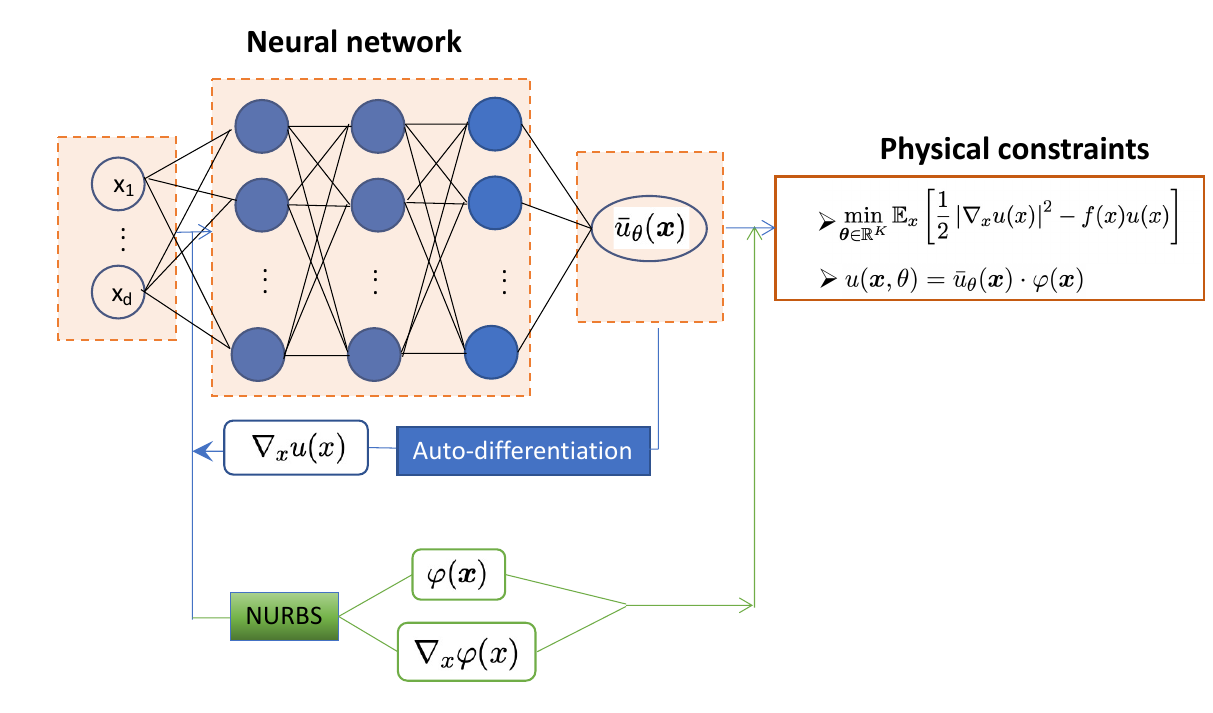}
    \caption{Deep NURBS method illustrated. Integration of a Neural Network with NURBS functions to generate admissible Neural Networks.}\label{fig:intro}
\end{figure}

One might argue that PINNs is a class of machine learning methods incorporating physical laws (PDEs) within their loss functions. Other traits characterizing different methods that obey the latter principle make them a variant of PINNs. Therefore, we believe the method we propose here can be considered a variant of PINNs in a manner that both are based upon the same principle of incorporating physical laws. The variational approach followed here relates more to the Deep Ritz method \cite{weinan2018} though.   

\subsubsection{Deep neural networks} Let us reconsider the physical domain $\cl{D}$. $\bs{x}$ is a random variable in that domain ($\bs{x} \in \cl{D}$) that can be sampled from a specific distribution that we will discuss in the sampling section. Whether we are willing to solve a PDE or to approximate an unknown function $f$, the neural network approximation proved to be an extremely efficient machine learning-based tool in either case \cite{hornik1989}. Denoting it $\bar{u}_{\theta}(\bs{x})$, the neural network approximation can take the following form:

\begin{equation}\label{eq:nn1}
\bar{u}_{\theta}(\bs{x})=\hat{\bs{\beta}} \cdot z_{\theta}(\bs{x})+b,
\end{equation}
where $\hat{\bs{\beta}}$ and $b$ denote the last layer NN parameters and $z_{\theta}(\bs{x})$ is a compositional function that has the following expression : 
\begin{equation}
z_{\theta}(\bs{x})=z_{l} \circ \ldots \circ z_{1}(\bs{x})
\end{equation}

$l$ is the number of hidden layers and $z_i$ is the output of the $i{th}$ layer that is given as function of its precedent layer's output, i.e $z_{i-1}$ as: $z_{i}=\sigma\left(\Theta_{i} \cdot z_{i-1}+b_{i}\right), \quad \Theta_{i} \in \mathbb{R}^{K \times K'}, b_{i} \in \mathbb{R}^{K}$
$\Theta$ is referred to as the weights of the NN, and $K$ ($K'$) is the size of the ith ($i-1$) layer (number of neurons). $\sigma$ is the non-linear function (activation function) that allows the approximating neural networks to cope with the complexity of the sought solution. A careful choice of the non-linear function has to take place for fast convergence and to prevent undesired effects related to the nature of the non-linear function. For example, Rectified linear units (ReLU) and Sigmoid functions that are known to be widely used activation functions suffer for deep architectures from the vanishing gradient issue \cite{Hochreiter1998}. Furthermore and from PINNs perspective, one has to choose an activation function with a certain degree of regularity that is adequate with the differentiation order of the PDE. We, hence, choose the function $\sigma = \max(0,x^3)$ that led to good results when applied to a similar framework \cite{weinan2018}.  We are interested in the forward scheme to solve the PDE assuming its well-definedness.To apply this pure unsupervised scheme, one does not require any labeled data from the domain $\cl{D}$. Unlike problems classified as forward problems, PINNs can refer to inverse or parameterized schemes in which the PDEs are equipped with uncertain parameters. The approach then uses labeled data to determine the best estimates for these parameters. 
\subsubsection{Inclusion of physical laws}
As known for NN - based approaches, a residual function needs to be defined so then its minimization leads to the true solution. The choice of this residual (loss) function is important to define the driving bias towards the best estimate. In supervised machine learning, most of the regression problems were defined using quadratic loss terms figuring the difference between the true data values and the NN predictions. That scheme proved useful for some applications but irrelevant to many science-related problems where the generalization (extrapolation) to out-of-the-range data is important. PINNs and their variants came with the idea of imposing physical constraints. For this purpose, one has to incorporate the residual form of the PDE given at Eq. \eqref{eq:pde} into the loss function.

\begin{equation}\label{eq:resd_abst}
\mathcal{R}\left(\bs{x} ; u; \frac{\partial u}{\partial x_{1}} \ldots; \frac{\partial^{2} u}{\partial x_{1}}, \ldots f \right)=0, \quad \bs{x} \in \cl{D}
\end{equation}

\noindent $\mathcal{R}$ is the residual of the PDE, which we chose to be the loss function for our optimization problem. $\mathcal{R}$ is given by the difference between the right and left-hand sides of the upper equation in Eq. \eqref{eq:pde}. This imposition of the PDE ensures the training process will be guided by the minimization of Eq. \eqref{eq:pde}. In other words, we are searching for the best trial function that fulfils that given PDE. This merge between the PDE residual and the loss function in ML terminology afforded PINNs unprecedented success in solving problems without data supervision.
\subsubsection{Deep NURBS}\label{sec:deepn}
Until this stage, nothing is done to ensure our problem's well-posedness (hence the solution's uniqueness). In the literature, boundaries have been treated differently depending on the approach used. In collocation methods, adding an extra term containing the boundary loss was more convenient, referred to as boundary operator (see the lower part of Eq. \eqref{eq:pde}). Approaches using losses in weak forms tend to add a parameterized integral to account for the boundaries (see Ref. \cite{weinan2018} in Eq. (3.2)). Shin et al. \cite{shin2020} reported that minimizers satisfying boundary conditions adopt a convergence mode $H^1$ compared to $C^0$ for minimizers that do not comply with boundary conditions. Accordingly, we define a method to incorporate boundary conditions inspired by (and based on) IsoGeometric Analysis. The method is known as Deep NURBS.\\
The overall idea of Deep NURBS is to develop our trial function as a product between an admissible field (constructed as in IsoGeometric Analysis \cite{Cot09}) and a NN estimate that is independent of boundary points: 
\begin{equation}
    u(\bs{x},\theta) = \varphi(\bs{x};\bs{p})\cdot \bar{u}(\bs{x},\theta)
\end{equation}
$\varphi(\bs{x};\bs{p})$ is the admissible field evaluated at the data points $\bs{x}$ for a given set of control points $\bs{p}$. We have intentionally omitted the $\bs{p}$ dependence of the total estimate as no optimization with respect to $\bs{p}$ is going to take place. Depending on the problem at hand, $\bs{p}$ can be defined from a simple guess. In the current case, the control points $\bs{p}$ are defined to concentrate sampled points at near complex regions of the domain. This will lead to an important feature, namely, an automatic importance sampling of the data.

\subsubsection{Sampling method (Importance sampling)}\label{sec:sampl_mthd}

To calculate the loss function given in Eq. \eqref{eq:weak_form} (Eq. \eqref{eq:resd_abst} for an abstract form), one has to search for an effective sampling procedure for a fast and convergent minimization.
Note that because we are using concepts from IsoGeometric Analysis, we use in the current approach two sampling methods. The first one occurs in the parametric space, which we refer to as parametric sampling, and physical space sampling is related to parametric sampling in a way that we will discuss in the next subsections.

\begin{figure}
    \centering
    \includegraphics[width=0.5\textwidth]{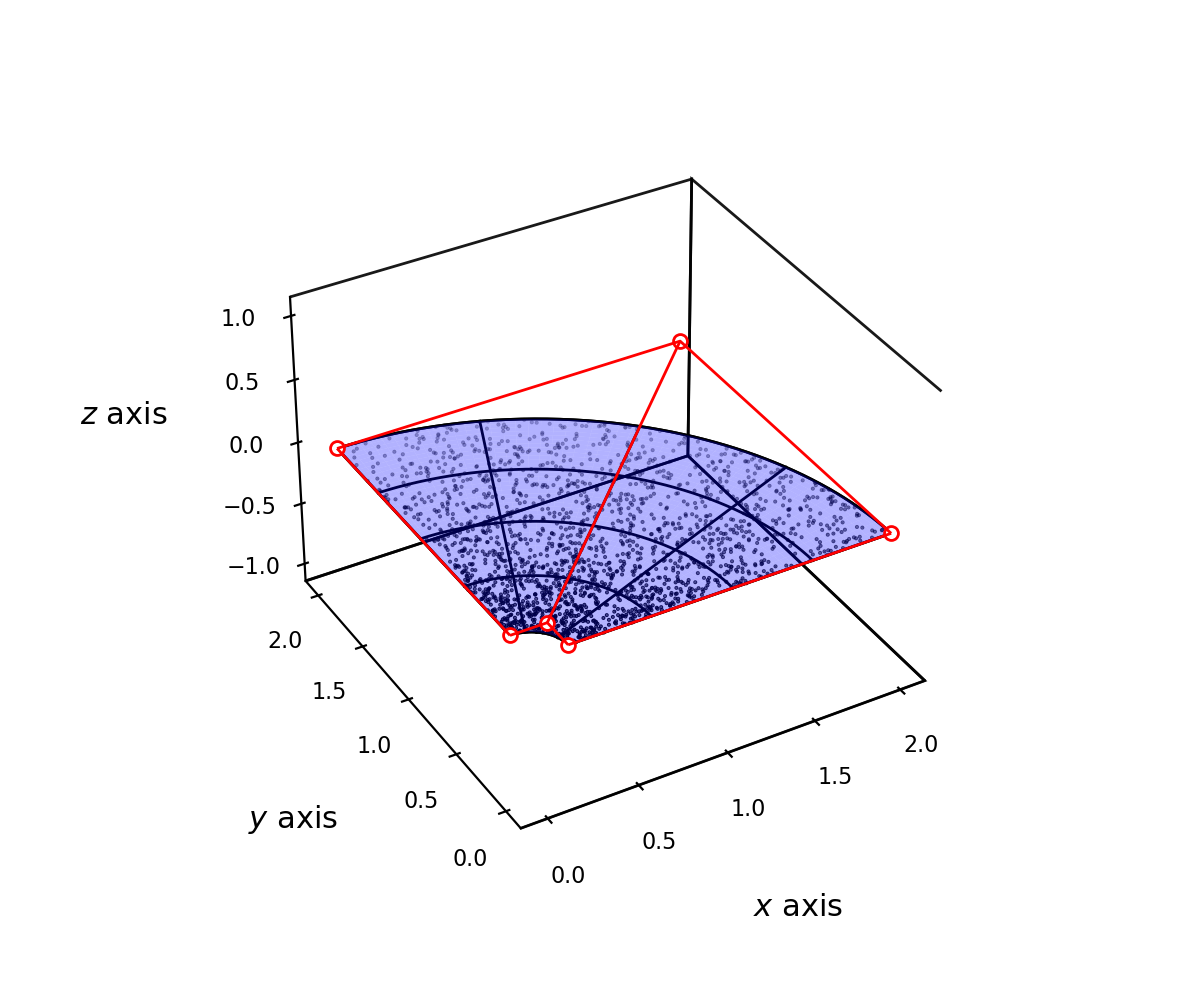}
    \caption{Parametric sampling.}
    \label{fg:parametric.sampling}
\end{figure}

\noindent\textit{Parametric sampling:} Motivated by quadrature rules for isoparametric elements in the finite element method (FEM), we here motivate our parametric sampling choice. In FEM, usually in regions with complex geometrical features, the geometry is discretized with sufficiently many elements. In each element, quantities of interest are integrated with a local Gaussian quadrature rule. Thus, in these regions, there is a concentration of quadrature points. In our current approach, we sample from the parametric space $\bs{\Xi}$ rather than $\cl{D}$, see Figure \ref{fg:parametric.sampling}.

\noindent\textit{Physical domain sampling:} The physical and parametric domains are related via the transformation $\chi$ given by : 
\begin{equation}\label{eq:transf}
\bs{\chi}(\bs{\xi})=\frac{1}{W(\boldsymbol{\xi})} \sum_{i=0}^{n_{1}} \sum_{j=0}^{n_{2}} N_{i, p_{1}}\left(\xi^{1}\right) N_{j, p_{2}}\left(\xi^{2}\right) \boldsymbol{p}_{i, j},
\end{equation}
This equation is equivalent to equation \eqref{eq:NURBS2d}; thus, with this transformation, it is possible to generate the physical domain given uniformly sampled data from the parametric space. Their respective distribution densities, i.e., the uniform parametric and the physical space distributions, are related via the Jacobian of the same transformation in the following manner \cite{asmussen2007}:
\begin{equation}\label{eq:density}
\rho_{\cl{D}}(\bs{x}) = \rho_{\bs{\Xi}}\left(\bs{\chi}^{-1}(\bs{x})\right)\left|\operatorname{det}\left[\left. \Grad_{\bs{y}}\bs{\chi}^{-1}(\bs{y})\right|_{\bs{y}=\bs{x}}\right]\right|.
\end{equation}

The argument of the determinant is the Jacobian of the inverse transformation $\bs{\chi}^{-1}$. $\rho_{\bs{\Xi}}(\bs{x})$ is the uniform parametric density distribution in the parametric domain and $\rho_{\cl{D}}(\bs{x})$ is the non-uniform density distribution in the physical domain. Looking at Eq. \eqref{eq:density}, one can see its dependence on the gradient of the transformation $\Grad_{\bs{y}}\bs{\chi}^{-1}(\bs{y})$ (since the transformation is bijective). The physical domain density will, therefore, depend on the concentration of the control points $\boldsymbol{p}_{i, j}$.  
Thus, by concentrating the control points at the boundary regions and close to singularities, one can have an automatic importance sampling in the physical domain. As can be seen in the Results section, the latter fact is behind the fast convergence of the Deep NURBS method.

Lastly, it is important to note that since we perform a parametric sampling and then map this to the physical domain, our Deep NURBS method can handle complex geometries simply and efficiently. Also, since our approximate is an admissible field, that is, Dirichlet boundary conditions are automatically satisfied and no sampling on the boundary is required. This technique drastically increases the approximability of Neural Networks.

\subsubsection{Optimization (Automatic differentiation)}\label{sec:optmz}

The evaluation of the energy (loss) term must go via the evaluation of the NN function's partial derivatives as a function of the input variables $\bs{x}$ after using the importance sampling described in the previous section. Fortunately, one can achieve this task in a rigorous way using automatic differentiation (AD) \cite{Rumelhart1986,wengert1964}. Furthermore AD is time-efficient in addition to being rigorous because it demands the evaluation of the function only twice, i.e., a forward and a backward evaluation of the target function $u$ regardless of the input variable dimension (see Ref. \cite{Lu2021} for an illustrative example). 

The PDE residual (our objective function) is expressed using derivatives derived using AD. An additional differentiation step is required for the training purpose, namely, differentiating the objective function as a function of the NN approximation parameters ($\theta$). We define these derivatives as $\nabla_{\theta}\cl{R}(\theta)$. The update rule of the training depends on the optimization order to take place. First-order optimizations such as the one used in this study have the following update rule \cite{adam2014}: 

\begin{equation}\label{eq:adam1}
\theta_{i+1}=\theta_i-\frac{\eta}{\sqrt{\hat{v}_i}+\epsilon} \hat{m}_i.
\end{equation}

$\epsilon$, $\beta_1$ and $\beta_2$ are hyperparamters where $\hat m_i$ and $\hat{v}_i$ denote the bias-corrected first and second moment estimates, respectively. These can be written as function of the error gradients as follows: 
\begin{equation}\label{eq:adam2}
\begin{aligned}
\hat{m}_i &=\frac{m_i}{1-\beta^i_1} = \frac{\beta_1 m_{i-1}+\left(1-\beta_1\right) \nabla_\theta \cl{R}(\theta_i)}{1-\beta^i_1} \\
\hat{v}_i &=\frac{v_i}{1-\beta^i_1} = \frac{\beta_2 v_{i-1}+\left(1-\beta_2\right) \left(\nabla_\theta \cl{R}(\theta_i)\right)^2}{{1-\beta^i_2}}
\end{aligned}
\end{equation}

$\theta_{i+1}$ in Eq. \eqref{eq:adam1} denotes the new variable affection after the $i$th iteration. The evaluation of $\cl{R}(\theta_i)$ is usually performed using mini-batches of data \cite{robbins1951}. $\gamma$ refers to the so-called learning rate. For the sake of convergence, the condition $\gamma_n \rightarrow 0$ as $n \rightarrow \infty$ needs to be satisfied \cite{asmussen2007}. A regularization term could have been appended to the residual figuring in Eq. \eqref{eq:adam2} in order to ensure the boundedness of the NN sought parameters.
Instead, the optimization is restricted to smaller search spaces, which can be seen as an implicit regularization \cite{shin2020}. In the next section, we will see the application of the Deep NURBS method to various problems with different types of boundaries, the results therein support the latter claim. 

\section{Numerical results}\label{sec:sec3}

In the previous section, we presented the overall framework of the deep NURBS (DN) method. In this section, we validate its efficiency for various cases with different geometries. The application of DN to real-world problems should be the focus of future work. We used three different examples with different boundary conditions and geometrical setups. As for the optimization, we used the default hyperparameters suggested in Adam's original paper \cite{adam2014}. The value of the learning rate $\gamma$ is set fixed to the value of $10^{-3}$ for all the examples. Additionally, all simulations were run with 300k epochs to guarantee that the comparisons are fair, that is, the optimization errors are negligible and the errors are due to the approximant characteristics. In the first example, we emphasize the manner in which the DN method copes with rough discontinuities and the presence of singularities. In the second example, we present the results obtained for a smoother geometry with a well-guessed admissible field $\varphi$ (Eq. \eqref{eq:admField}). In the last example, we focus on the impact of the choice of $\varphi$ on the optimization convergence rate and on the estimate accuracy. The latter fact we refer to as the inductive bias.
We consider the two-dimensional Poisson equation given in Eq. \eqref{eq:poisson1}

\subsection{Physical Domain with Corner Singularity }\label{sec:exple1}

\begin{figure}
    \centering
    \includegraphics[width=0.5\textwidth]{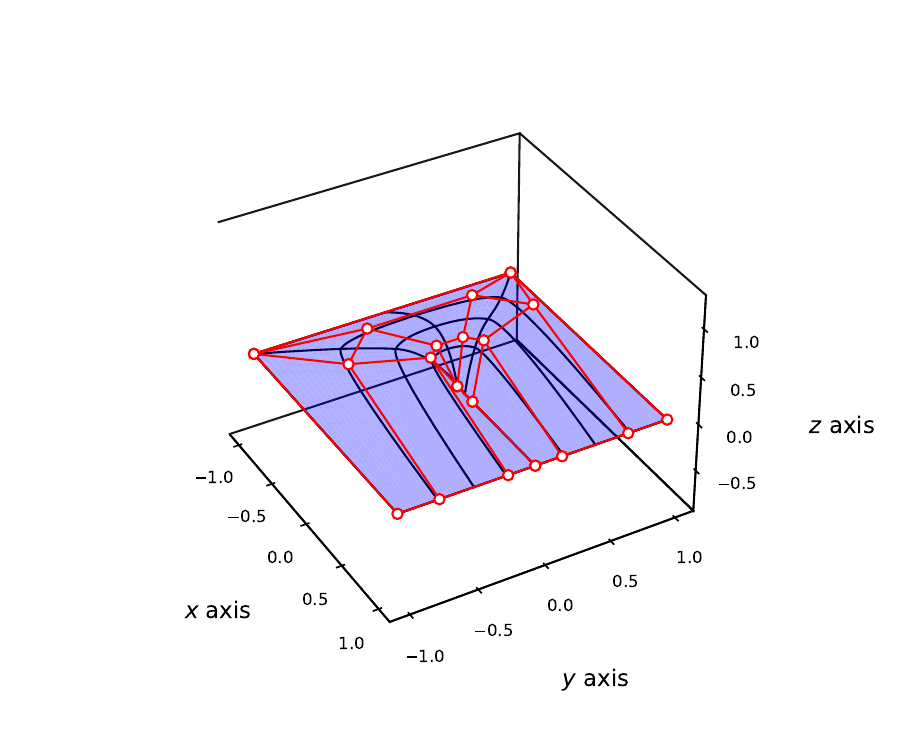}
    \caption{NURBS parameterization for the homogeneous Dirichlet boundary condition in $\cl{D}$ for problem \ref{sec:exple1}. NURBS: degree $p_1=p_2=2$, knot vectors: $\Xi^1=\{0,0,0,0.25,0.5,0.75,1,1,1\}$ and $\Xi^2=\{0,0,0,0.5,1,1,1\}$.}
    \label{fig:wein_arch}
\end{figure}

\begin{figure}[ht!]\label{weinan4p}
    \centering
\includegraphics[width=.7\textwidth]{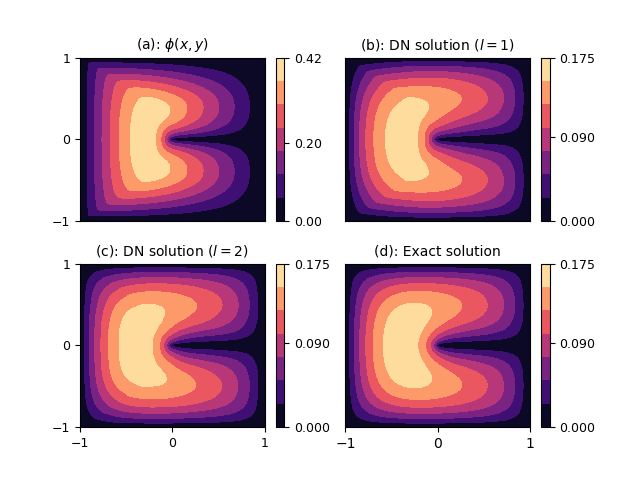}
    \caption{Solution of the Poisson equation for problem \ref{sec:exple1}. We use the admissible field $(\phi,\bs{\zeta}=\bs{0})\in(\cl{P},\cl{Z})$ (a) to predict the solution using one (b) and two (c) hidden layers Deep NURBS solver. These results are to be compared to the exact solution (d) calculated using FEM}
    \label{fig:weinan4p}
\end{figure}

In this example, we are interested in solving the Poisson's equation defined as follows:

\begin{equation}\label{eq:poisson2}
\begin{aligned}
-\Delta u(\bs{x})=1, & \quad \bs{x} \in \Omega \\
        u(\bs{x})=0, & \quad \bs{x} \in \partial \Omega,\\
\end{aligned}
\end{equation}
where the physical domain is given by $\cl{D}=(-1,1) \times(-1,1) \backslash[0,1) \times\{0\}$. $\Omega$ exhibits a sharp discontinuity along the segment $[0,1)$ and singularity at the origin. The NURBS geometry is presented in Fig. \ref{fig:wein_arch}. The control points (referred to as $\bs{p}_{i,j}$ in Eq. \eqref{eq:transf}) are selected at the boundaries and the neighborhood of the origin to impose boundary conditions and to prompt importance sampling. The admissible field in Fig.\ref{fig:weinan4p}(a) has arbitrary values except at the boundaries, i.e., $\varphi$ is chosen to impose the Dirichlet boundary conditions. The smoothness of $\varphi$ inside the domain $\Omega$ is controlled based on the degree of the NURBS functions (as per Eq. \eqref{eq:admField}). The admissible field's regularity impacts the minimizers' convergence to the PDE solution. Note that the mapping is not invertible at the origin; however, this point is never sampled as it lies on the boundary.\\
\noindent We address the current problem using two different NN architectures separately and independently (only for comparison) with one and two hidden layers, respectively. Note that the number of NN - layers used throughout all the paper examples does not exceed two.  $\max(0,x^3)$ is used as the main activation function for the hidden layers, whereas the linear activation function is used for the output layer. None of the NN regularizing techniques (such as dropout) are used here. To ensure that the method is generalizable and not sensitive to parameter adjustments, we intentionally implement the basic version of multilayer perceptrons.

\noindent Indeed, a good level of accuracy can be achieved due to the implementation of physical constraints and the hard imposition of the boundary conditions using NURBS. Moreover, the automatic importance sampling which is enhanced by DN (see sec.\ref{sec:optmz}) plays a critical role in promoting fast convergence. The number of neurons per layer is fixed to the value of $50$. The loss function (PDE residual) is given by:

$$\cl{R} = \mathbb{E}_{\bs{x}}\left[\frac{1}{2}\left|\nabla_{\bs{x}} u(\bs{x})\right|^{2}- u(\bs{x})\right]$$

$\bs{x}=\left(x_1,x_2\right)$ is a 2D variable. Using importance sampling, the expectation of the latter equation can be transformed into a weighted sum: 
\begin{equation}
    \cl{R} \approx \frac{1}{N}\sum^{N}_{i=1}\frac{\left|\nabla_{\bs{x}} u(\bs{x}_i)\right|^{2}- u(\bs{x}_i)}{\rho_{\mathcal{D}}(\bs{x}_i)}
\end{equation}
$N$ is the number of samples, and $\rho_{\mathcal{D}}(\bs{x})$ is the physical domain probability density. $\rho_{\mathcal{D}}(\bs{x})$ is related to the transformation $\bs{\chi}$, which maps the parametric space samples to the ones in the physical domain. The transformation $\bs{\chi}$ is given in Eq. \eqref{eq:transf}, whereas $\rho_{\mathcal{D}}(\bs{x})$ is calculated using Eq. \eqref{eq:density}. We use first-order stochastic optimization implemented in Adam optimizer \cite{adam2014}.

\begin{table}[ht!]
\caption{Errors for problem \ref{sec:exple1} using the Deep Ritz and the Deep NURBS methods (after 300k epochs). Both methods use a fully connected NN architecture composed of two hidden layers containing 50 neurons each.}\label{tab:tab1}
\centering
\begin{tabular}{lccc}
\hline Approach & MSE & Relative $L_{2}$ & $L_{\infty}$ \\
\hline Deep NURBS & $2.8 \times 10^{-5}$ & $2.5 \times 10^{-3}$ & $2.3 \times 10^{-2}$ \\
Deep Ritz & $3 \times 10^{-3}$ & $8.7 \times 10^{-1}$ & $1.1 \times 10^{-1}$
\end{tabular}
\end{table}

\noindent Solutions to the PDE of Eq. \eqref{eq:poisson2} for the geometry illustrated in Fig. \ref{fig:wein_arch} are shown in Fig. \ref{fig:weinan4p}. We present the solutions to the PDE for two different NNs with one and two hidden layers, respectively, alongside the exact solution. Table. \ref{tab:tab1} show a comparative study between DN and classical Deep Ritz method \cite{weinan2018} (DR) performances. Both methods used two hidden layers in their NN approximations. We used the hyperbolic tangent (tanh) activation function for DR, except that we used the same shareable parameters for both methods. \\

\begin{figure}[ht!]
    \centering
    \includegraphics[width=.95\textwidth]{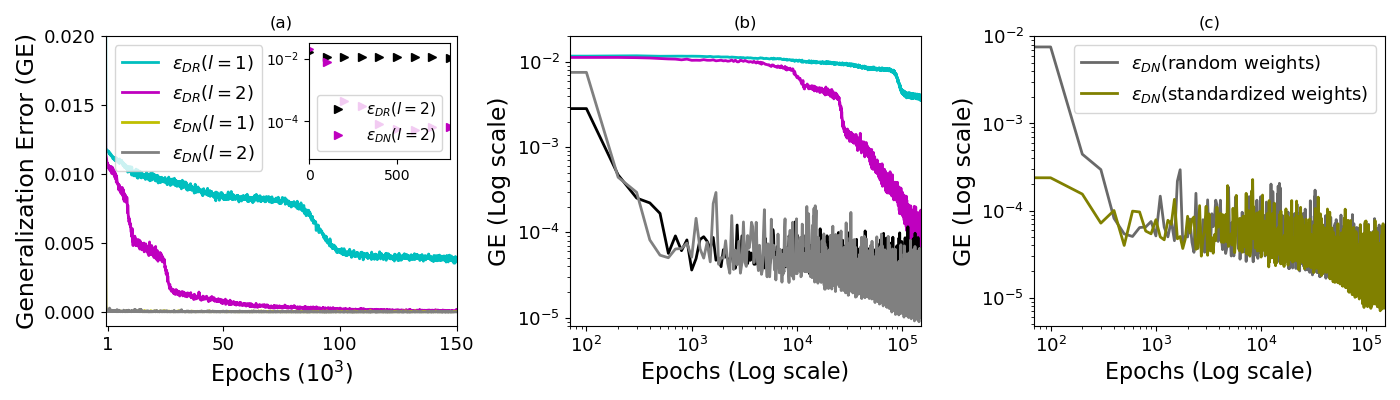}
    \caption{Training errors for Poisson's equation solution to problem \ref{sec:exple1}. (a) shows the error lines obtained using Deep NURBS and Deep Ritz methods for one and two hidden layers. In the inset, the logarithmic scale is applied to the y-axis for the case of neural networks with two hidden layers. (b) is the all axes log-scaled version of (a). In (c), we plotted the error lines for our method with random initial weights (gray) and with weights pre-trained to the identity estimator (olive).}\label{fig:train_history}
\end{figure}

\noindent The DN method outperforms DR by two orders of magnitude for the Mean Squared Error (MSE) and the relative $L_2$ error, and by one order of magnitude for the $L_\infty$ error. This dissimilarity is because, with shallow NN, DR could not penalize the energy term of the left-hand side of Eq. \eqref{eq:weak_form} with the boundary term in the right-hand side of the same equation. Moreover, the natural imposition of boundary conditions within the DN framework resulted in a truncation of the search space to a smaller space (\cite{shin2020}). Because the complexity of the NN approximation is related to the number of its layers and neurons, one concludes that the solution can be, in this case, approximated by fewer layers of NNs. The DN approximation for one or two hidden layers has similar error ranges. DN, with only one hidden layer, seems to exhibit slightly sharp edges originating from the original admissible field. A choice of a smooth admissible field can lead to more accurate results. Note that by \textit{smooth} here, we mean the choice of a higher polynomial degree for the admissible field shown in Eq. \eqref{eq:admField}.

\begin{figure}[ht!]
    \centering
    \includegraphics[width=.7\linewidth]{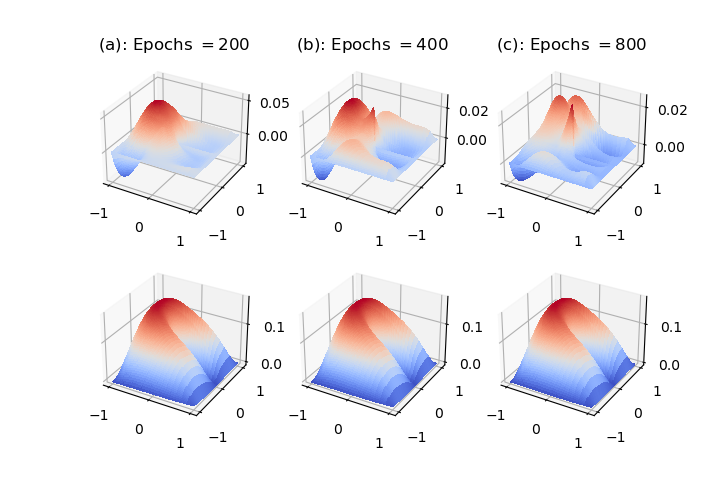}
    \caption{Three-dimensional error plots calculated during the learning process for DN (upper curves) and DR (lower curves) for different numbers of epochs.}
    \label{fig:error_plots_wein}
\end{figure}

\subsubsection{Impact of variance reduction} DN qualitative convergence rate can be inherently deduced from Fig.\ref{fig:train_history}((a) and (b) subplots). The importance sampling being a variance reduction technique has led to faster convergence of the DN method compared to classical DR. Inset of Fig.\ref{fig:train_history}(a) shows that an error of less than $10^{-4}$ is obtained using DN compared to a stagnated value around $10^{-2}$ of the MSE metric in the case of DR. An error value as less as $10^{-4}$ is obtained in less than $500$ epochs. As a comparative study, we show in (c) the error lines for our method estimations of the solution using random initial (gray curve) and pre-trained (olive-colored curve) weights. In the latter case, the NN parameters are obtained by training a separate NN to predict the identity output. Such a choice will ensure that the initial guess would have the same symmetry as the admissible field.

Note that we selected this example due to its complex (sharply discontinued) geometry. This particular example has been addressed \cite{weinan2018} using the Deep Ritz method with a quite involved residual network architecture. Using our method, we could solve the same PDE with only one (or two) hidden layers NN. The latter insight is supported by the error contours shown in Fig. \ref{fig:error_plots_wein}.\\

\subsection{Annulus geometry}\label{sec:annl_geo}

\begin{figure}[ht!]
    \centering
    \includegraphics[width=0.5\textwidth]{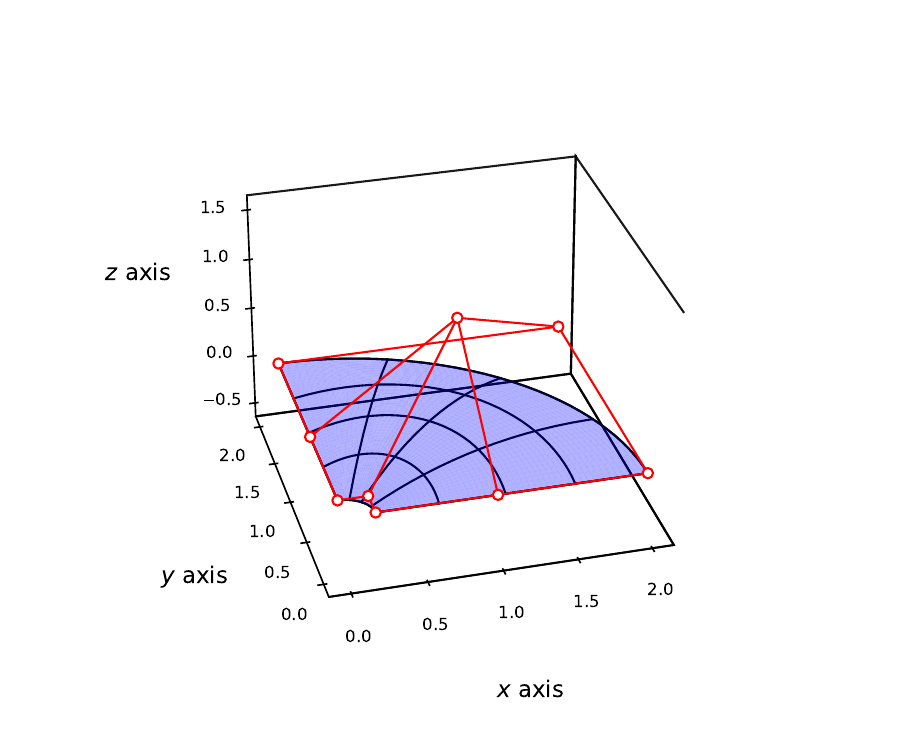}
    \caption{NURBS parameterization for homogeneous Dirichlet boundary condition in $\cl{D}$ for the problem \ref{sec:annl_geo}. NURBS: degree $p_1=p_2=2$, knot vectors: $\Xi^1=\{0,0,0,1,1,1\}$ and $\Xi^2=\{0,0,0,1,1,1\}$.}
    \label{fig:annulus}
\end{figure}

The geometry in this example consists of a 2D annulus quarter. In polar coordinates, the physical domain is spanned by the vector ($r$,$\phi$) where $r \in \left[r_1,r_2\right]$  and $\phi \in [0,\frac{\pi}{2}]$. We selected $r_1$ and $r_2$ to be equal to $0.2$ and $2$, respectively. The boundary conditions over the region $\delta\Omega$ are given as: 

$$
\begin{aligned}
u_{\theta}\left(r, \phi\right) &=0,  \quad  r=\{r_1,r_2\}, \phi \in[0,\frac{\pi}{2}] \\
u_{\theta}\left(r, \phi\right) &=0, \quad r=[r_1,r_2],  \phi \in \{0,\frac{\pi}{2}\} \\
\end{aligned}
$$

\begin{figure}[ht!]
    \centering
    \includegraphics[width=0.65\textwidth]{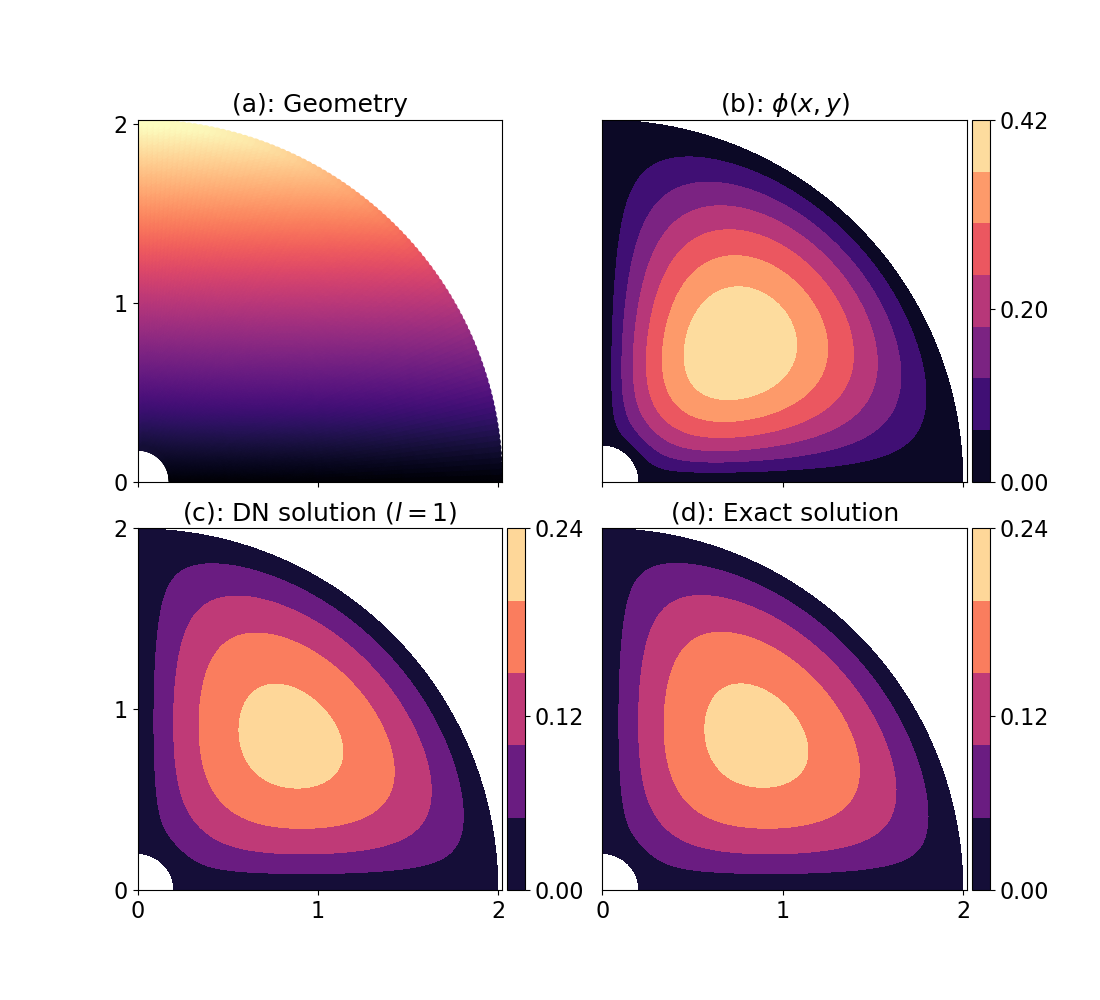}
    \caption{Solution of the Poisson equation for the annulus problem for which the geometry is shown in (a). (b) showcases the corresponding admissible field $(\phi,\bs{\zeta}=\bs{0})\in(\cl{P},\cl{Z})$. The Deep NURBS solution (c) is calculated using one hidden layer with 50 neurons. The reference (exact) solution is given in (d)}
    \label{fig:4plots_ann}
\end{figure}

We choose to solve the same Poisson equation stated in Eq. \eqref{eq:poisson2}. The NURBS geometry is shown in Fig.\ref{fig:annulus}. The control points $\bs{p}_{i,j}$ are placed along the boundary regions and inside the physical domain. The solution to the PDE alongside the chosen admissible field is given in Fig. \ref{fig:4plots_ann}. Being a multiplicative factor, the admissible field is expected to induce an inductive bias towards the values of its generating function $\varphi$. 
\noindent The NN approximation $\bar{u}(\bs{x},\theta)$ tends, during the training, to oppose this bias given the PDE and the boundary values. Table. \ref{tab:tab2} shows different error metrics for the DN approximation of the sought solution. The error values are calculated in reference to an exact FEM calculation. An MSE error of the order of $10^{-7}$ is reported for the DN estimator with one and two hidden layers. Moreover, the increase of the NN complexity ($l=2$) did not alter the generalization error. With an extremely small relative $L_2$ error of the order of $10^{-6}$, one can confirm that the estimation accuracy is no more dependent on the NN complexity. This makes the current method a serious preference for higher dimensional problems.          

\begin{table}[ht!]
\caption{Errors for problem \ref{sec:annl_geo} (after 300k epochs) using the Deep NURBS method for one and two hidden layers. The method uses a fully connected NN architecture with 50 neurons per layer.}\label{tab:tab2}
\centering
\begin{tabular}{lclc}
\hline Number of layers & $\mathrm{MSE}$ & Relative $L_{2}$ & $L_{\infty}$ \\
\hline $\mathrm{l}=1$ & $1.1 \times 10^{-7}$ & $8.4 \times 10^{-6}$ & $9.3 \times 10^{-4}$ \\
$\mathrm{l}=2$ & $1.2 \times 10^{-7}$ & $8.8 \times 10^{-6}$ & $9.1 \times 10^{-4}$
\end{tabular}
\end{table}

\subsection{Square with a hole}\label{sec:exple3}
\begin{figure}[ht!]
    \centering
    \includegraphics[width=0.5\textwidth]{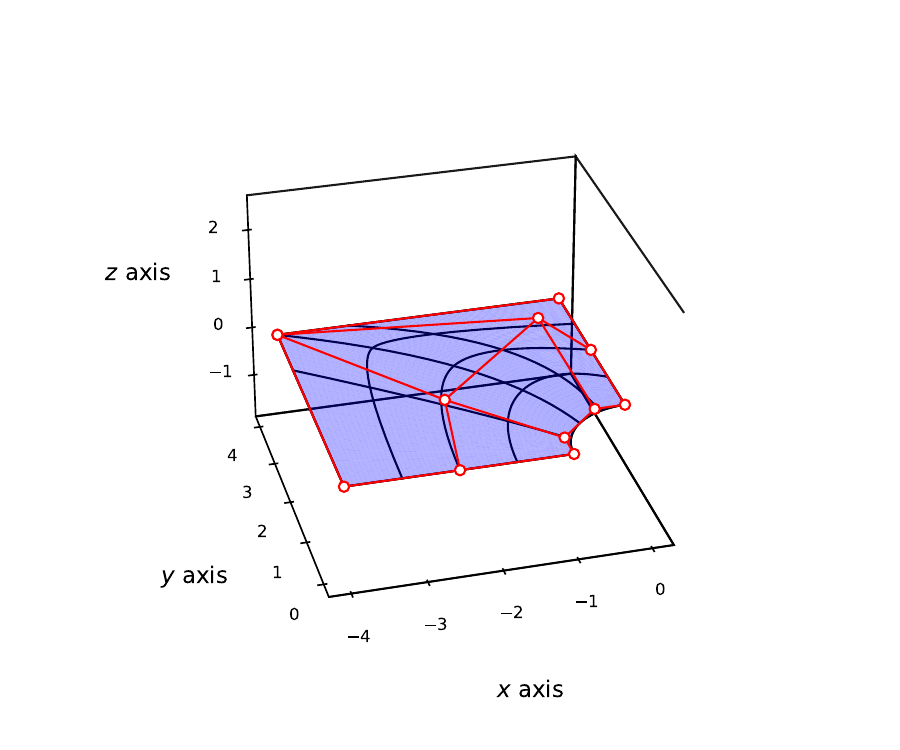}
    \caption{NURBS parameterization for homogeneous Dirichlet boundary condition in $\cl{D}$ for the problem  \ref{sec:exple3}. NURBS: degree $p_1=p_2=2$, knot vectors: $\Xi^1=\{0,0,0,0.5,1,1,1\}$ and $\Xi^2=\{0,0,0,1,1,1\}$.}
    \label{fig:square}
\end{figure}

\begin{figure}[ht!]
    \centering
    \includegraphics[width=.5\textwidth]{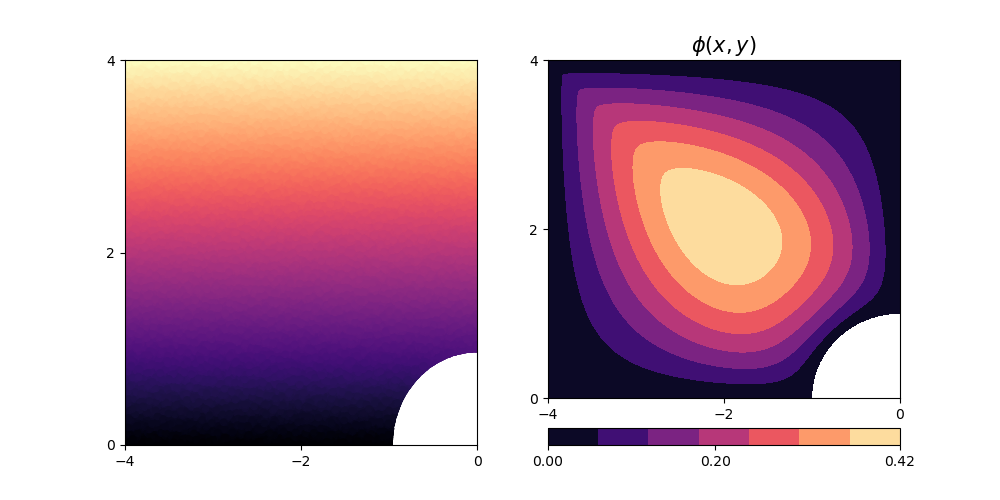}
    \caption{(a) shows the physical domain for the problem \ref{sec:exple3}, where (b) shows the corresponding admissible field $(\phi,\bs{\zeta}=\bs{0})\in(\cl{P},\cl{Z})$ }
    \label{fig:quarter_arch}
\end{figure}

This example comprises a square disc (with length $x_d = 4$) with a hole (with a radius $r = 1$) as shown in Fig. \ref{fig:quarter_arch}. The NURBS geometry is shown in Fig.\ref{fig:square}.

\noindent In addition to the boundary conditions imposed by the disc shape, an additional discontinuity is brought out by the circular hole along its perimeter. Therefore the admissible field $\varphi$ obeys the following Dirichlet boundary conditions: 
$$
\begin{aligned}
\varphi(x,y) &= 0 \quad (x,y) = \{-4,4\}\cup \{4,-4\}, \\
\varphi(x,y) &= 0 \quad x^2 + y^2 = 1
\end{aligned} $$

\begin{figure}[ht!]
    \centering
    \includegraphics[width=.9\textwidth]{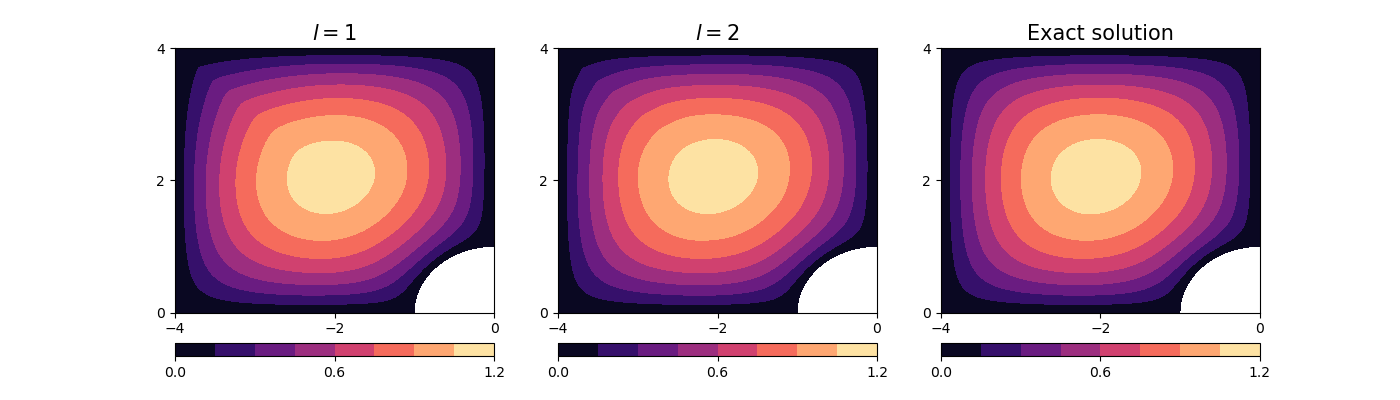}
    \caption{Solution of the Poisson equation for the problem \ref{sec:exple3} using one (a) and two (b) hidden layers NN solver. The reference solution is showcased in (c).}
    \label{fig:3plots_quart}
\end{figure}

\noindent The admissible field plot is shown in Fig. (\ref{fig:quarter_arch} - (b)). Despite complying with the problem boundary conditions, the admissible field $\varphi(\bs{x})$ of Fig.\ref{fig:quarter_arch} - (b) exhibits a noticeable directional anisotropy that induces an inductive bias affecting the training cost and efficiency of the NN multiplying function $\bar{u}(\bs{x},\theta)$.
This can be revealed from a comparison between the admissible field plot in Fig.\ref{fig:quarter_arch} and the solution for one hidden layer NN plotted in Fig.\ref{fig:3plots_quart} - (a). A twofold solution is required to be implemented to cope with this situation. The choice of the control point locations should depend on a prior understanding of the PDE symmetry and geometry. Moreover, one can mitigate the impact of the inductive bias by selecting more regular admissible fields. Using higher order polynomials (by increasing the number of B-splines basis functions in Eq. \eqref{eq:N}) produces fields that are indeed not hard to adjust by the multiplying NN function $\bar{u}(\bs{x},\theta)$.

\begin{figure}[ht!]
    \centering
    \includegraphics[width=0.9\textwidth]{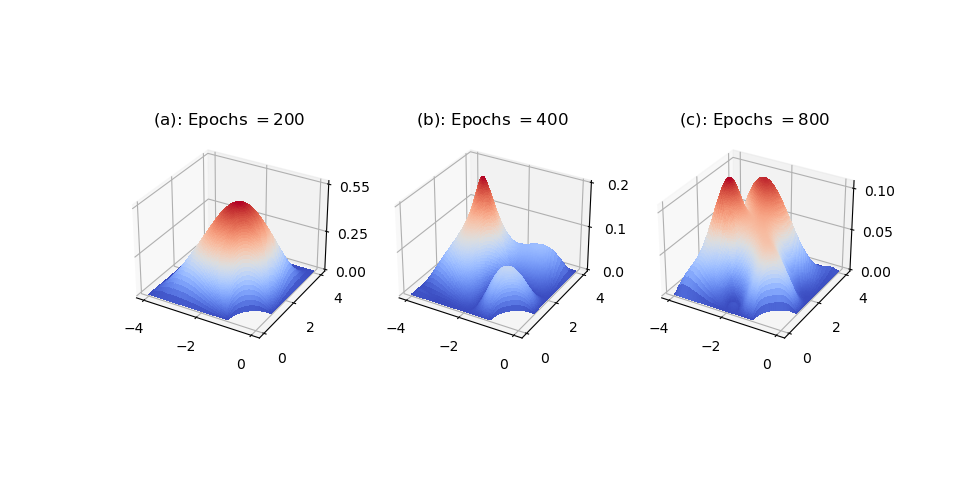}
    \caption{Three-dimensional training error plots obtained at different epochs for the solution of the PDE in Eq. \eqref{eq:poisson2} for the physical domain presented in Fig. \ref{fig:quarter_arch}}
    \label{fig:err_quart}
\end{figure}

Despite the biased admissible field, the NN has been able to cope with this discrepancy with a shallow NN (l=1 and l=2), as can be seen in Fig.\ref{fig:3plots_quart}. Three-dimensional error plots can be showed from Fig. \ref{fig:err_quart}. 

\begin{table}[ht!]
\caption{Errors for problem \ref{sec:exple3} (after 300k epochs) using the Deep NURBS method for one and two hidden layers. The method uses a fully connected NN architecture with 50 neurons per layer.}\label{tab:tab3}
\centering
\begin{tabular}{lclc}
\hline Number of layers & $\mathrm{MSE}$ & Relative $L_{2}$ & $L_{\infty}$ \\
\hline $\mathrm{l}=1$ & $2.9 \times 10^{-5}$ & $7.2 \times 10^{-5}$ & $2.1 \times 10^{-2}$ \\
$\mathrm{l}=2$ & $2.5 \times 10^{-5}$ & $6.2 \times 10^{-6}$ & $8.7 \times 10^{-3}$
\end{tabular}
\end{table}

Table \ref{tab:tab3} shows the approximation errors for the deep NURBS method for a NN with one and two hidden layers. Unlike the previous example, the increase in the number of layers in this case ($l=2$) resulted in a decrease of the generalization error by one order of magnitude. This is primarily because of the inductive bias induced by the admissible field of Fig. \ref{fig:quarter_arch}. A careful choice of the admissible field can be, nonetheless, performed effortlessly by figuring out factors that might affect the solution, including but not limited to the problem symmetry. \\
To summarize, the three examples showed the applicability and efficiency of the Deep NURBS algorithm in solving PDE forward problems. We have seen that the method can interestingly solve PDEs with harsh boundary conditions like the one in the first example (sec. \ref{sec:exple1}) with the minimum NN setup. DN proves to be a fast and convergent method for moderately complex boundary conditions, a fact that is revealed in the second example. Example 3 shows how DN can cope with \textit{wrongly biased} priors. Nevertheless, a good choice of priors would positively impact its capability and efficiency.    

\section{Conclusion}

We proposed in this study a new method based on combining PINNs and IsoGeometric analysis. The method (which we refer to as Deep NURBS) combines the simple implementation of PINNs and its capability to solve PDEs with the power of the NURBS method for imposing Dirichlet boundary conditions even for complex geometries. We have described how the use of DN can result in an automatic importance sampling-induced variance reduction. We have demonstrated that the method can be used to solve PDEs on complex geometries with extremely shallow NNs by applying it to various problems with various geometries. In comparison to traditional PINNs methods, the automatic variance reduction had an impact on the optimization, resulting in effective training with fewer iterations. Note that we have used our method for problems having different geometries as we believe this new method is meant to deal with complex geometries, a common feature of many real problems. The application of DN to more concrete and real-world problems is the scope for future work.

\section{Declarations: Funding and/or Conflicts of interests/Competing interests \& Data availability}

This publication is based upon work supported by the King Abdullah University of Science and Technology (KAUST) Office of Sponsored Research (OSR) under Award No. OSR-2019-CRG8-4033 and the Alexander von Humboldt Foundation.

The authors have no conflicts to disclose. Codes can be obtained by reasonable request to the corresponding author.

\section{Additional Declaration}

All simulations were ran in an iMac 3,8 GHz 8-core Intel Core i7, 128 GB 2667 MHz DDR4.





\footnotesize


\begin{thebibliography}{10}

\bibitem{LeCun2015}
Y~LeCun, Y~Bengio, and G~Hinton.
\newblock Deep learning.
\newblock {\em Nature}, 521:436--444, 2015.

\bibitem{speech2012}
G~Hinton, L~Deng, D~Yu, GE~Dahl, A-R Mohamed, N~Jaitly, A~Senior, V~Vanhoucke, P~Nguyen, TN~Sainath, and B~Kingsbury.
\newblock Deep neural networks for acoustic modeling in speech recognition: The shared views of four research groups.
\newblock {\em IEEE Signal Processing Magazine}, 29(6):82--97, 2012.

\bibitem{speech2013}
A~Graves, A-R Mohamed, and G~Hinton.
\newblock Speech recognition with deep recurrent neural networks.
\newblock In {\em 2013 IEEE International Conference on Acoustics, Speech and Signal Processing}, pages 6645--6649, 2013.

\bibitem{img_recg2012}
A~Krizhevsky, I~Sutskever, and G~Hinton.
\newblock Imagenet classification with deep convolutional neural networks.
\newblock {\em Proc. Advances in Neural Information Processing Systems}, 25, 2012.

\bibitem{img_recg2015}
C~Szegedy, W~Liu, Y~Jia, P~Sermanet, S~Reed, D~Anguelov, D~Erhan, V~Vanhoucke, and A~Rabinovich.
\newblock Going deeper with convolutions.
\newblock In {\em 2015 IEEE Conference on Computer Vision and Pattern Recognition (CVPR)}, pages 1--9, 2015.

\bibitem{lagaris1998}
IE~Lagaris, A~Likas, and DI~Fotiadis.
\newblock Artificial neural networks for solving ordinary and partial differential equations.
\newblock {\em IEEE Transactions on Neural Networks}, 9(5):987--1000, 1998.

\bibitem{carleo2017}
G~Carleo and M~Troyer.
\newblock Solving the quantum many-body problem with artificial neural networks.
\newblock {\em Science}, 355:602, February 2017.

\bibitem{carleo2018}
G~Carleo, Y~Nomura, and M~Imada.
\newblock Constructing exact representations of quantum many-body systems with deep neural networks.
\newblock {\em Nature Communications}, 9:5322, 2018.

\bibitem{Torlai2018}
G~Torlai, G~Mazzola, J~Carrasquilla, M~Troyer, R~Melko, and G~Carleo.
\newblock Neural-network quantum state tomography.
\newblock {\em Nature Physics}, 14:447--450, 2018.

\bibitem{schmidt2009}
M~Schmidt and H~Lipson.
\newblock Distilling free-form natural laws from experimental data.
\newblock {\em Science}, 324(5923):81--85, 2009.

\bibitem{snyder2012}
J~Snyder, M~Rupp, K~Hansen, K~M\"uller, and K~Burke.
\newblock Finding density functionals with machine learning.
\newblock {\em Phys. Rev. Lett.}, 108:253002, Jun 2012.

\bibitem{raissi2019}
M~Raissi, P~Perdikaris, and GE~Karniadakis.
\newblock Physics-informed neural networks: A deep learning framework for solving forward and inverse problems involving nonlinear partial differential equations.
\newblock {\em Journal of Computational Physics}, 378:686--707, 2019.

\bibitem{berg2018}
J~Berg and K~Nyström.
\newblock A unified deep artificial neural network approach to partial differential equations in complex geometries.
\newblock {\em Neurocomputing}, 317:28--41, 2018.

\bibitem{weinan2018}
E~Weinan and Y~Bing.
\newblock The deep ritz method: A deep learning-based numerical algorithm for solving variational problems.
\newblock {\em Communications in Mathematics and Statistics}, 6:1--12, 2018.

\bibitem{sirignano2018}
J~Sirignano and K~Spiliopoulos.
\newblock Dgm: A deep learning algorithm for solving partial differential equations.
\newblock {\em Journal of Computational Physics}, 375:1339--1364, 2018.

\bibitem{cai2021}
S~Cai, Z~Wang, L~Lu, TA~Zaki, and GE~Karniadakis.
\newblock Deepm\&mnet: Inferring the electroconvection multiphysics fields based on operator approximation by neural networks.
\newblock {\em Journal of Computational Physics}, 436:110296, 2021.

\bibitem{Kurth2018}
T~Kurth, S~Treichler, J~Romero, M~Mudigonda, N~Luehr, E~Phillips, A~Mahesh, M~Matheson, J~Deslippe, M~Fatica, P~Prabhat, and M~Houston.
\newblock Exascale deep learning for climate analytics.
\newblock In {\em SC18: International Conference for High Performance Computing, Networking, Storage and Analysis}, pages 649--660, 2018.

\bibitem{samaniego2020}
E~Samaniego, C~Anitescu, S~Goswami, VM~Nguyen-Thanh, H~Guo, K~Hamdia, X~Zhuang, and T~Rabczuk.
\newblock An energy approach to the solution of partial differential equations in computational mechanics via machine learning: Concepts, implementation and applications.
\newblock {\em Computer Methods in Applied Mechanics and Engineering}, 362:112790, 2020.

\bibitem{Han2018}
J~Han, A~Jentzen, and E~Weinan.
\newblock Solving high-dimensional partial differential equations using deep learning.
\newblock {\em Proceedings of the National Academy of Sciences}, 115(34):8505--8510, 2018.

\bibitem{Lu2021}
L~Lu, P~Jin, G~Pang, Z~Zhang, and GE~Karniadakis.
\newblock Learning nonlinear operators via deeponet based on the universal approximation theorem of operators.
\newblock {\em Nature Machine Intelligence}, 3(3):218--229, 2021.

\bibitem{Bru20}
SL~Brunton, BR~Noack, and P~Koumoutsakos.
\newblock Machine learning for fluid mechanics.
\newblock {\em Annual Review of Fluid Mechanics}, 52:477--508, 2020.

\bibitem{mishra2022}
S~Mishra and R~Molinaro.
\newblock {Estimates on the generalization error of physics-informed neural networks for approximating PDEs}.
\newblock {\em IMA Journal of Numerical Analysis}, 2022.

\bibitem{fang2021}
Z~Fang.
\newblock A high-efficient hybrid physics-informed neural networks based on convolutional neural network.
\newblock {\em IEEE Transactions on Neural Networks and Learning Systems}, pages 1--13, 2021.

\bibitem{shin2020}
Y~Shin.
\newblock On the convergence of physics informed neural networks for linear second-order elliptic and parabolic type pdes.
\newblock {\em arXiv: Numerical Analysis}, 2020.

\bibitem{sukumar2022}
N~Sukumar and A~Srivastava.
\newblock Exact imposition of boundary conditions with distance functions in physics-informed deep neural networks.
\newblock {\em Computer Methods in Applied Mechanics and Engineering}, 389:114333, 2022.

\bibitem{tensorflow2016}
M~Abadi, P~Barham, J~Chen, Z~Chen, A~Davis, J~Dean, M~Devin, S~Ghemawat, G~Irving, M~Isard, et~al.
\newblock Tensorflow: A system for large-scale machine learning.
\newblock In {\em 12th $\{$USENIX$\}$ Symposium on Operating Systems Design and Implementation ($\{$OSDI$\}$ 16)}, pages 265--283, 2016.

\bibitem{Deb72}
Carl De~Boor.
\newblock On calculating with b-splines.
\newblock {\em Journal of Approximation theory}, 6(1):50--62, 1972.

\bibitem{Cox72}
Maurice~G Cox.
\newblock The numerical evaluation of b-splines.
\newblock {\em IMA Journal of Applied mathematics}, 10(2):134--149, 1972.

\bibitem{Bez72}
Pierre B{\'e}zier.
\newblock Numerical control-mathematics and applications.
\newblock {\em Translated by AR Forrest}, 1972.

\bibitem{hornik1989}
Kurt Hornik, Maxwell Stinchcombe, and Halbert White.
\newblock Multilayer feedforward networks are universal approximators.
\newblock {\em Neural Networks}, 2(5):359--366, 1989.

\bibitem{Hochreiter1998}
Sepp Hochreiter.
\newblock The vanishing gradient problem during learning recurrent neural nets and problem solutions.
\newblock {\em International Journal of Uncertainty, Fuzziness and Knowledge-Based Systems}, 06(02):107--116, 1998.

\bibitem{Cot09}
J~Austin Cottrell, Thomas~JR Hughes, and Yuri Bazilevs.
\newblock {\em Isogeometric analysis: toward integration of CAD and FEA}.
\newblock John Wiley \& Sons, 2009.

\bibitem{asmussen2007}
Søren Asmussen and Peter~W. Glynn.
\newblock Stochastic simulation : algorithms and analysis.
\newblock {\em Springer}, 2007.

\bibitem{Rumelhart1986}
David~E. Rumelhart, Geoffrey~E. Hinton, and Ronald~J. Williams.
\newblock Learning representations by back-propagating errors.
\newblock {\em Nature}, 323(6088):533--536, 1986.

\bibitem{wengert1964}
R.~E. Wengert.
\newblock A simple automatic derivative evaluation program.
\newblock {\em Commun. ACM}, 7(8):463–464, aug 1964.

\bibitem{adam2014}
Diederik~P Kingma and Jimmy Ba.
\newblock Adam: A method for stochastic optimization.
\newblock {\em arXiv preprint arXiv:1412.6980}, 2014.

\bibitem{robbins1951}
Herbert Robbins and Sutton Monro.
\newblock {A Stochastic Approximation Method}.
\newblock {\em The Annals of Mathematical Statistics}, 22(3):400 -- 407, 1951.

\end{thebibliography}

\end{document}